\documentclass{amsart}
\let\<=\langle
\let\>=\rangle
\let\noi=\noindent
\let\sse=\subseteq
\let\vphi=\varphi
\let\veps=\varepsilon
\let\limply=\Longrightarrow

\def\0{\{0\}}
\def\span{{\kern.5pt{\rm span}\kern1pt}}
\def\smallfrac#1#2{{\textstyle{\frac{#1}{#2}}}}

\def\conv{{\;\longrightarrow\;}}
\def\wconv{{{\buildrel_{\scriptstyle w}\over\conv}}}
\def\sconv{{{\buildrel_{\scriptstyle s}\over\conv}}}
\def\uconv{{{\buildrel_{\scriptstyle u}\over\conv}}}
\def\notwconv{{{\wconv\kern-13pt\sslash}\kern9pt}}
\def\notsconv{{{\sconv\kern-13pt\sslash}\kern9pt}}
\def\notuconv{{{\uconv\kern-13pt\sslash}\kern9pt}}

\font\fiverm=cmr5
\def\sslash{\hbox{{\fiverm /}}}
\def\Ps{{\vbox{\hbox{$\wp\kern.3pt$}\vskip.2pt}}}

\def\B{{\mathcal B}}
\def\H{{\mathcal H}}

\def\T{{\mathcal T}}
\def\X{{\mathcal X}}

\def\BX{{\B[\X]}}
\def\BH{{\B[\H]}}

\def\CC{{\mathbb C\kern.5pt}}
\def\DD{{\mathbb D\kern.5pt}}
\def\JJ{{\mathbb J\kern.5pt}}
\def\QQ{{\mathbb Q\kern.5pt}}
\def\NN{{\mathbb N\kern.5pt}}
\def\RR{{\mathbb R\kern.5pt}}
\def\TT{{\mathbb T\kern.5pt}}
\def\ZZ{{\mathbb Z\kern.5pt}}

\newsymbol\varnothing 203F
\let\void=\varnothing

\def\smallmatrix#1{\null\,\vcenter{
                   \baselineskip=8pt\mathsurround=0pt\ialign{
                   \hfil ${\scriptstyle##}$
                   \hfil &&
                   \hfil ${\scriptstyle##}$
                   \hfil \crcr
                   \mathstrut \crcr
                   \noalign{\kern-\baselineskip}#1 \crcr
                   \mathstrut \crcr
                   \noalign{\kern-\baselineskip} \crcr }}\!}

\begin{document}

\vglue-72pt\noi
\hfill{\it Advances in Mathematical Sciences and Applications}\/
{\bf 35}(2) (2026) 559--578

\vglue15pt
\title
      [Weak Quasistability and Rajchman measures]
      {Weak Quasistability and Rajchman measures}
\author
       [C.S. Kubrusly]
       {Carlos Kubrusly}
\address{Catholic University of Rio de Janeiro, Brasil}
\email{carlos@ele.puc-rio.br}
\renewcommand{\keywordsname}{Keywords}
\dedicatory{Dedicated to the memory Ren\'ee Good\/ $($1988--2026\/$)$
            and Alex Pretti\/ $($1988--2026\/$)$}
\date{Received September 1, 2025; Accepted January 29, 2026}
	
\begin{abstract}
It is shown that weak quasistability does not imply power boundedness, but
coercive power unbounded operators cannot be weakly quasistable.\ Although a
finite measure over the unit disc is a Rajchman measure if and only if the
position operator is weakly stable, it is shown that the position operator is
weakly quasistable for every finite continuous measure over the unit disc.\
Corollaries linking Rajchman measures with weak stability and weak
quasistability follow the above results.\
\end{abstract}

\maketitle

\vskip-10pt\noi
\section{Introduction}
Stability of an operator $T$ (e.g., acting on a normed space) means
convergence of the power sequences $\{T^n\}$ to the null operator.\ Weak
stability holds when this happens in the weak sense.\ Weak quasistability
holds when convergence is weakened to convergence of a subsequence.\ For
recent expositions on weak stability and weak quasistability, see \cite{JJKS}
and \cite{KV}, where these notions are discussed in detail and their
importance is explained and exemplified.\ The present paper investigates the
relation of weak stability and weak quasistability with Rajchman measures.\

\vskip6pt
The main results proved here are:
\vskip6pt
\begin{description}
\item{$\kern-2pt\circ$}
Weakly quasistable operators are either power bounded or noncoercively power
unbounded (Theorem 4.3).\
\vskip4pt
\item{$\kern-2pt\circ$}
If $\{z^k\}$ is an orthonormal basis for $L^2(\TT,\mu)$, then the measure
$\mu$ is \hbox{Rajchman}, and the converse fails:\ there are Rajchman measures
$\mu$ for which any pair of distinct elements from $\{z^k\}$ is not orthogonal
in $L^2(\TT,\mu)$ (Theorem 6.1).\
\vskip4pt
\item{$\kern-2pt\circ$}
The position operator on $L^2(\TT,\mu)$ is weakly quasistable for every
continuous measure $\mu$ (Theorem 7.3).\
\end{description}
\vskip6pt\noi
Among other conclusions, we highlight the following ones:\
\begin{description}
\vskip6pt
\item{$\kern-2pt\circ$}
Weakly quasistability does not imply power boundedness (Remark~4.4).\
\vskip4pt
\item{$\kern-2pt\circ$}
The position operator on $L^2(\TT,\mu)$ is weakly stable if and only if
$\mu$ is a \hbox{Rajchman} measure, and therefore weak stability is ensured
by testing the unit function only (Proposition 5.3).\
\end{description}

\vskip6pt
The paper is organised into 6 more sections.\ Basic notation and terminology
are outlined in Section 2.\ The definitions of weak stability and weak
quasistability are posed and compared in Section 3.\ Section 4 analyses when
a weakly quasistable operator is power bounded.\ The relation of Rajchman
measures with weak \hbox{stability} of the position operator, which is
unitary, is explored in Section 5.\ Section 6 \hbox{discusses} the connection
of Rajchman measures $\mu$ and orthonormal bases for ${L^2(\TT,\mu)}.$
\hbox{Section} 7 closes the paper by showing that every finite continuous
measure over the unit circle is quasi-Rajchman, which implies weak
quasistability for the position operator, and also characterises non-Rajchman
continuous measures.\

\section{Preliminaries:\ Notation and Terminology}

The following common notation will be adopted throughout the paper.\ The set
of all integers, the set of all positive integers, the real line, and the
complex plane are denoted by $\ZZ$, $\NN$, $\RR$, and $\CC$, respectively, as
usual.\ Let $\DD$, $\DD^-\!$, and $\TT$ stand for the open unit disc, the
closed unit disc, and the unit circle in $\CC$, respectively.\

\vskip6pt
A sequence of positive integers is any $\NN$-valued function on $\NN.$ Regard
the set of all positive integers equipped with its natural well-ordering as
the self-indexed sequence of all positive integers so that
$\NN=\{n\}_{n\ge1}.$ A subsequence $\{n_j\}=\{n_j\}_{j\ge1}$ of the positive
integers is a strictly increasing (infinite) sequence of positive integers.\

\vskip6pt
We will identify a subsequence of the positive integers (i.e., a strictly
increasing function of $\NN$ into itself) with its range.\ So a subsequence of
the positive integers is identified with an infinite, ordered, strictly
increasing subset of $\NN.$ We will use the same notation, $\{n_j\}$, for
either of them.\

\vskip6pt
A subsequence of the positive integers is {\it nontrivial}\/ if there is
another (equally infinite) subsequence such that they have no common entries.\
If $\{n'_j\}$ is a subsequence of a subsequence $\{n_j\}$ of the positive
integers, then $\{n_j\}$ is a {\it supersequence}\/~of~$\{n'_j\}.$

\vskip6pt
A subsequence $\{a_{n_j}\}$ of an $A$-valued sequence $\{a_n\}$ (for an
arbitrary nonempty set $A$) is the restriction of $\{a_n\}$ to a subsequence
$\{n_j\}$ of the positive integers.\

\vskip6pt
A subsequence $\{n_j\}$ of the positive integers is of {\it bounded
increments}\/ (or {\it bounded gaps}\/) if ${\sup_j(n_{j+1}\!-n_j)<\infty}.$ A
subsequence $\{a_{n_j}\}$ of $\{a_n\}$ is {\it boundedly spaced}\/ if it is
indexed by a subsequence $\{n_j\}$ of bounded increments.\

\vskip6pt
Let $\X$ be a normed space and let $\X^*\!$ be its dual.\ An $\X$-valued
sequence $\{x_n\}$ is weakly convergent if there exists an ${x\in\X}$ such
that ${\lim_nf(x_n)=f(x)}$ for \hbox{every} ${f\in\X^*}$; equivalently, such
that ${\lim_nf(x_n\!-x)=0}$ for every ${f\in\X^*\kern-1pt}.$ In this case the
sequence $\{x_n\}$ is said to converge weakly to the vector $x.$ An
alternative and usual notation for weak convergence that will be used here is
${x_n\!\wconv x}\,$ or $\,x\kern-1pt=\kern-1pt{w\hbox{\,-}\lim_n x_n}.$

\vskip6pt
The above definition is standard.\ Actually, ``every normed space has a
topology $\T$ [the weak topology on it] such that a sequence in the space
converges weakly to an element of the space if and only if the sequence
converges to that element with respect to $\T.$ For the moment, the statement
that a sequence converges weakly to a certain limit should not be taken to
imply anything more than is stated [above]'' \cite[p.116]{Meg}.\ We will not
deal with weak topology techniques here (although the above definition
coincides with convergence in the weak topology in a Hilbert space).\

\section{Weak Quasistability:\ Definition and Comparison}

By an operator on a normed space $\X$ we mean a bounded linear (i.e., a
continuous linear) transformation of $\X$ into itself.\ Let $\BX$ stand for
the normed algebra of~all operators on $\X.$ The same notation ${\|\;\;\|}$
will be used for the norm on $\X$ and for the induced uniform norm on $\BX.$
An operator $T$ is {\it strongly stable}\/ if ${\lim_n\!\|T^nx\|=0}$ for every
${x\in\X}$ (notation: $\!{T^n\!\sconv O}$), and {\it uniformly stable}\/ if
${\lim_n\!\|T^n\|=0}$ (nota\-tion: $\!{T^n\!\uconv O}).$ An operator $T$ is
{\it power bounded}\/ if ${\sup_n\!\|T^n\|\kern-1pt<\!\infty}$; equivalently,
if ${\sup_n\!\|T^nx\|\kern-1pt<\!\infty}$ for every $x$ in $\X$ by the
Banach--Steinhaus Theorem if $\X$ is a \hbox{Banach} space.\ If $\X$ is a
\hbox{Hilbert space}, then it will be denoted by $\H$, and ${T^*\!\in\BH}$
will stand for the adjoint of ${T\kern-1pt\in\BH}.$

\vskip6pt\noi
{\bf Definition 3.1$.$}
(a)
An operator ${T\kern-1pt\in\BX}$ is {\it weakly stable}\/ if, for each
${x\in\X}$,
$$
{\lim}_n|f(T^nx)|=0
\;\;\hbox{for every}\;\,
f\in\X^*\!.
$$
Notation: ${w\hbox{\,-}\lim_n T^n x=0}$, or ${T^nx\wconv 0}$, for every
${x\in\X};\,$ or simply $\;{T^n\!\wconv O}.$
\vskip6pt\noi
(b)
Equivalently, an operator $T$ on a normed space $\X$ is {\it weakly stable}\/
if, for each vector ${x\in\X}$ and every subsequence $\{n_j\}$ of the
positive integers,
$$
{\lim}_j|f(T^{n_j}x)|=0
\;\;\hbox{for every}\;\,
f\in\X^*\!.
$$
(c)
A subsequence $\{n_j\}$ of the positive integers is a {\it subsequence of
weak stability}\/ of $T$ for $x$ if ${\lim_j|f(T^{n_j}x)|=0}$ for every $f.$
So a weakly stable operator is one for which every subsequence of the positive
integers is of weak stability for every~vector~$x.$

\vskip6pt
An operator is weakly quasistable if the limit in Definition 3.1(a) is
weakened to limit inferior.

\vskip6pt\noi
{\bf Definition 3.2$.$}
(a)
An operator ${T\kern-1pt\in\BX}$ is {\it weakly quasistable}\/ if, for each
${x\in\X}$,
$$
{\liminf}_n|f(T^nx)|=0
\;\;\hbox{for every}\;\,
f\in\X^*\!.
$$
Notation: ${w\hbox{\,-}\liminf_nT^nx=0}$ for every ${x\in\X}.$
\vskip6pt\noi
(b)
Equivalently, an operator $T$ on a normed space $\X$ is {\it weakly
quasistable}\/ if for each vector ${x\in\X}$ there exists a subsequence
$\{n_j\}=\{n_j(x)\}$ of the positive integers (that may depend on each $x$
but does not depend on $f$) for which
$$
{\lim}_j|f(T^{n_j}x)|=0
\;\;\hbox{for every}\;\,
f\in\X^*\!.
$$
(c)
So a weakly quasistable operator is an operator $T$ such that, for every $x$,
there is at least one subsequence ${\{n_j\}}$ (that may depend on $x$) of weak
stability of $T$~for~$x.$

\vskip6pt
Weak quasistability plays an important role in weak l-sequential supercyclicity
\cite[Corollaries 4.3 and 4.4]{KD}.\

\vskip6pt
It is worth noticing that the notion of quasistability has a plain meaning
only in the weak case.\ In fact, it is easy to see (cf.\
\cite[Propositions 4.1 and~4.2]{KV}$\kern.5pt$) that
$$
\hbox{$\limsup_n\|T^n\|=0\iff\lim_n\|T^n\|=0$}
$$
(i.e., quasistability always coincides with stability for the uniform case)
and
$$
\hbox{
$\!$if $\sup_n\kern-1pt\|T^n\|\kern-1pt<\kern-1pt\infty$, then
$\limsup_n\kern-1pt\|T^nx\|\kern-1pt=\kern-1pt0\;\,\forall\,x\in\kern-1pt\X
\iff\,\lim_n\kern-1pt\|T^nx\|\kern-1pt=\kern-1pt0\;\,\forall\,x\in\kern-1pt\X$
}
$$
(i.e., for power-bounded operators, quasistability coincides with stability
for the strong case).\ For the weak case the notions are different even for
power-bounded operators \cite[Proposition 4.3]{KV}.\ However, if a weakly
quasistable operator $T$ is such that for every $x$ there is a subsequence
$\{n_j\}$ of weak stability that is of bounded increments, then $T$ is weakly
stable \cite[Theorem~5.3]{KV}.\

\vskip6pt\noi
{\bf Definition 3.3$.$}
(a)
An operator ${T\kern-1pt\in\BX}$ is {\it homogeneously weakly quasistable}\/
if it is weakly quasistable and there exists a subsequence $\{n_j\}$ of the
positive integers such that, for all vectors ${x\in\X}$,
$$
{\lim}_j|f(T^{n_j}x)|=0
\;\;\hbox{for every}\;\,
f\in\X^*\!.
$$
(b)
Equivalently, there is a subsequence of weak stability common to all
${x\in\X}.$
\vskip6pt\noi
(c)
So a homogeneously weakly quasistable is a weakly quasistable operator $T$
such that at least one subsequence ${\{n_j\}}$ of weak stability of $T$ does
not depend~on~$x.$

\vskip6pt
It is clear that every weakly stable operator is homogeneously
weakly quasistable, and every homogeneously weakly quasistable operator is
weakly quasistable.\ An~example of a power-bounded operator that is not weakly
stable but is homogeneously weakly quasistable will be presented in
Remark 7.4.

\vskip6pt\noi
{\bf Definition 3.4$.$}
(a)
An operator $T$ on a normed space $\X$ is {\it weakly unstable}\/ if
there exists an ${x_0\in\X}$ such that
$$
|f_0(T^nx_0)|\not\to0
\;\;\hbox{as}\;\:
n\to\infty
\;\;\hbox{for some}\;\,
f_0\in\X^*\!.
$$
Notation: ${T^n\!\notwconv O}$.
\vskip6pt\noi
(b)
Equivalently, ${T\kern-1pt\in\BX}$ is {\it weakly unstable}\/ if there exists
a vector ${x_0\in\X}$, a functional ${f_0\in\X^*\!}$, and a subsequence
$\{m_j\}$ of the positive integers
such that
$$
|f_0(T^{m_j}x_0)|\not\to0
\;\;\hbox{as}\;\;
j\to\infty.
$$
(c)
A subsequence $\{m_j\}$ of the positive integers is a {\it subsequence of weak
instability}\/ of $T$ for some vector $x_0$ if it is not a subsequence of weak
stability of $T$ for $x_0.$ A weakly unstable operator is one for which there
exists at least one subsequence of the positive integers that is of weak
instability for some vector $x_0.$

\vskip6pt\noi
{\sc Remark 3.5$.$}
{\sf The particular case of operators acting on a Hilbert space I.}
\vskip6pt\noi
If $\{T_n\}$ is a $\BX$-valued sequence of operators on a complex inner
product space $\X$ with inner product ${\<\;\,\,;\;\,\>}$, then the
polarisation identity ensures that
\begin{eqnarray*}
\<T_nx\,;y\>
&\kern-6pt=\kern-6pt&
\smallfrac{1}{4}\big(\<T_n(x+y)\,;(x+y)\>-\<T_n(x-y)\,;(x-y)\>             \\
&\kern-6pt+\kern-6pt&
i\<T_n(x+i\kern.5pty)\,;(x+i\kern.5pty)\>
-\<T_n(x-i\kern.5pty)\,;(x-i\kern.5pty)\>\big)
\end{eqnarray*}
for every ${x,y\in\X}$, which implies that the following assertions are
equivalent.
\vskip2pt\noi
\begin{description}
\item{$\kern-18pt$(a)$\kern4pt$}
${\lim_n|\<T_nx\,;y\>|=0}\;$ for every $\,{x,y\in\X}$.
\vskip6pt
\item{$\kern-18pt$(b)$\kern4pt$}
${\lim_n|\<T_nx\,;x\>|=0}\;$ for every $\,{x\in\X}$.
\end{description}
\vskip2pt\noi
uch an equivalence still holds in a real Hilbert space if each $T_n$ is
self-adjoint.\

\vskip6pt\noi
{\sc Remark 3.6$.$}
{\sf The particular case of operators acting on a Hilbert space II.}
\vskip6pt\noi
The Riesz Representation~Theorem in a Hilbert space ensures that if $\H$ is a
Hilbert space with inner product ${\<\;\,\,;\;\,\>}$, then according to
Definitions 3.1 and 3.2, an operator $T$ on $\H$ is weakly stable if and only
if 
$$
\kern-2pt\hbox{for each $x\kern-1pt\in\kern-.5pt\H$,
$\lim_n\kern-1pt|\<T^nx\,;y\>|\kern-1pt=\kern-1pt0$ for every
$y\kern-1pt\in\kern-.5pt\H$
\kern6pt(i.e., for every ${x,y\in\kern-.5pt\H}$).\kern-6pt}    \leqno{\rm(a)}
$$
As we saw in Remark 3.5, if $\H$ is a complex Hilbert space, then this is
equivalent~to
$$
\hbox{$\lim_n|\<T^nx\,;x\>|=0$ for every $x\in\H$.}            \leqno{\rm(b)}
$$
Similarly, an operator $T$ on a Hilbert $\H$ is weakly quasistable if and only
if
$$
\kern-6pt\hbox{for each $x\kern-1pt\in\kern-.5pt\H$,
$\liminf_n\kern-1pt|\<T^nx\,;y\>|\kern-1pt=\kern-1pt0$ for every
$y\kern-1pt\in\kern-.5pt\H$
\kern4pt(i.e., for every ${x,y\kern-1pt\in\kern-.5pt\H}$).\kern-12pt}
                                                               \leqno{\rm(c)}
$$
That is, for each ${x\in\H}$, there is a subsequence $\{n_j\}$ of the positive
integers such~that
$$
\hbox{$\lim_j|\<T^{n_j}x\,;y\>|=0$ for every $y\in\H$.}        \leqno{\rm(d)}
$$
\vskip1pt

\vskip6pt\noi
{\sc Remark 3.7$.$}
{\sf The particular case of operators acting on a Hilbert space III.}
\vskip6pt\noi
(a)
Weak stability for an operator $T$ on a Hilbert space $\H$ means that the
power sequence $\{T^n\}$ satisfies the equivalent limiting conditions (a) or
(b) in Remark 3.6.\ Strong stability for an operator $T$ on a Hilbert space
$\H$ means ${\lim_n\|T^nx\|\!=0}$ for every ${x\in\H}$, where ${\|\;\;\|}$ is
the norm generated by the inner product ${\<\;\,\,;\;\,\>}$, which clearly
implies (but is not implied by) weak stability.\ Differently from the strong
stability case, an operator $T$ and its adjoint $T^*$ are weakly stable
together.
\vskip5pt\noi
(b)
For weak quasistability we cannot expect in general a simplified counterpart
with ${y\kern-1pt=\kern-.5ptx}$, as in the weak stability case of
Remark 3.6(b).\ The reason~is~this$:$ for an arbitrary weakly quasistable
operator $T\kern-1pt$, the subsequence $\{n_j\}$ appearing in Remark 3.6(d)
may depend on $x.$ Indeed, the polarisation identity displayed in Remark 3.5
only ensures that for each ${x\in\H}$ there is a subsequence
$\{n_j(x)\}$~such~that
\begin{eqnarray*}
\<T^{n_j(x)}x\,;y\>
&\kern-6pt=\kern-6pt&
\smallfrac{1}{4}\big(\<T^{n_j(x)}(x+y)\,;(x+y)\>
-\<T^{n_j(x)}(x)(x-y)\,;(x-y)\>                                            \\
&\kern-6pt+\kern-6pt&
i\<T^{n_j(x)}(x+i\kern.5pty)\,;(x+i\kern.5pty)\>
-\<T^{n_j(x)}(x)(x-i\kern.5pty)\,;(x-i\kern.5pty)\>\big)
\end{eqnarray*}
for every ${y\in\H}$, and the condition ${\lim_n\<T^{n_j(x)}x\,;x\>=0}$ for
every ${x\in\H}$ does not guarantee that the right-hand side goes to zero for
every ${y\in\H}.$ For homogeneous weak stability, however, a simplified
counterpart with ${y\kern-1pt=\kern-.5ptx}$, as in the weak stability case of
Remark 3.6(b), holds naturally.\

\section{Weak Quasistability and Power Boundedness}

The following examples set a starting point for a discussion on power
boundedness and weak quasistability.

\vskip5pt\noi
{\sc Example 4.1$.$}
{\sf A power-bounded, weakly unstable but weakly quasistable operator.}

\vskip5pt\noi
Consider the Foguel operator
${F\kern-1pt=\kern-1pt\big(\smallmatrix{S^* & P               \cr
                                        O   & S^{\phantom{*}} \cr}\big)}$
on the direct sum ${\H\oplus\H}$ of a separable infinite-dimensional Hilbert
space $\H$ with itself.\ Here $S$ is a unilateral shift of multiplicity one
acting on $\H$ that shifts an orthonormal basis $\{{\bf e}_k\}_{k\ge0}$ for
$\H.$ The operator ${P\!:\H\!\to\H}$ is the orthogonal projection onto the
closure of the span of $\{{\bf e}_j\!:j\in\JJ\}$, where $\JJ$ is any sparse
infinite subset of positive integers with the following property: if
${i,j\in\JJ}$ and ${i<j}$, then ${2i<j}$ (e.g.,
$\JJ={\{j\ge 1\!:\;j=3^k;\;k\ge 0\}}$, the set of all integral powers of 3, is
a sample of a sparse set of positive integers satisfying the above
property).\ The operator $F$ was the first example of a power-bounded
operator that is not similar to a contraction \cite{Fog,Hal1}.\ Moreover, it
is also known that $F$ is not weakly stable.\ Since, for every ${n\ge0}$,
$$
F^n=\big(\smallmatrix{S^{*n} & P_n              \cr
                      O      & S^{n\phantom{*}} \cr}\big)
\quad\hbox{with}\quad
P_{n+1}={\sum}^n_{i=0}S^{*\,n-i}PS^i
\;\;\hbox{and}\;\;
P_0=O,
$$
\vskip-3pt\noi
we get
\vskip-3pt\noi
$$
\<F^{n+1}x\,;y\>=\<S^{*n+1}x_1,\;y_1\>+\<P_{n+1}x_2\,;y_1\>+\<S^{n+1}x_2,y_2\>
$$
\vskip3pt\noi
for every ${x=(x_1,x_2)\in\H\oplus\H}$ and every ${y=(y_1,y_2)\in\H\oplus\H}.$
As $S$ and $S^*$ are both weakly stable, it follows that
$$
{\liminf}_n|\<F^nx\,;y\>|=0
\quad\iff\quad
{\liminf}_n|\<P_nx_2\,;y_1\>|=0.
$$
It was shown in \cite[Proposition 4.3]{KV} that the power-bounded and weakly
unstable operator $F$ is weakly quasistable by exhibiting a subsequence
$\{n_j\}$ of the positive integers such that each $P_{n_j}x_2$ is orthogonal
to $y_1$ for every ${x_2,y_1\kern-1pt\in\H}$, and so
$\{\<P_{n_j}x_2\,;y_1\>\}$ is a null subsequence.\

\vskip5pt\noi
{\sc Remark 4.2$.$}
{\sf Spectral properties and power boundedness in weak quasistability.}

\vskip5pt\noi
(a)
As in the weak stability case, weak quasistability is also preserved under
the adjoint operation.\ An operator $T\kern-1pt$ and its adjoint $T^*\!$,
both acting on a Hilbert space $\H$, are weakly quasistable together because
$$
{\liminf}_n|\<T^nx\,;y\>|={\liminf}_n|\<T^{*n}y\,;x\>|
\;\;\hbox{for every}\;\,
x,y\in\H.
$$

\vskip0pt\noi
(b)
Let $\{\sigma_{\kern-1ptP}(T),\sigma_{\kern-1ptR}(T),\sigma_{\kern-.5ptC}(T)\}$
be the classical partition of the spectrum $\sigma(T)$ of\/ $T\kern-1pt$
consisting of point spectrum, residual spectrum, and continuous spectrum,
respectively.\ Since $T\kern-1pt$ and $T^*\!$ are weakly quasistable together
and
${\sigma_{\kern-1ptR}(T)
=\sigma_{\kern-1ptP}(T^*)^*\backslash\sigma_{\kern-1ptP}(T)}$, with
$\Lambda^*\!$ denoting the set of all complex conjugates of an arbitrary set
${\Lambda\kern-1pt\subset\kern-1pt\CC}$, it follows that
${\liminf_n\kern-1pt|\<T^nx\,;y\>|\kern-1pt=\kern-1pt0}$ for every
${x,y\kern-1pt\in\kern-1pt\H}$ implies
${\big(\sigma_R(T)\cup\sigma_{\kern-1ptP}(T)\kern.5pt\big)\cap\TT
\kern-1pt=\kern-1pt\void}.$
Thus, as in the case of weak stability,
$$
\hbox{$T$ is weakly quasistable $\;\limply\;$
$\sigma(T)\cap\TT\sse\sigma_{\kern-.5pt C}(T)$.}
$$
The Gelfand--Beurling formula for the spectral radius shows at once that a
power-bounded operator is such that ${\sigma(T)\sse\DD^-\!}.$ Therefore,
$$
\hbox{$T$ is weakly quasistable and power bounded $\;\limply\;$
$\sigma(T)\cap\TT\sse\sigma_{\kern-.5pt C}(T)\sse\DD^-\!$.}
$$
\vskip-2pt

\vskip6pt
$\!$If an operator $T\kern-1pt$ on a Banach space is power unbounded (i.e.,
not power bounded), then there is an ${x_0\kern-1pt\in\kern-1pt\X}$ for which
${\sup_n\|T^nx_0\|\kern-1pt=\kern-1pt\infty}$ by the Banach--Steinhaus
Theorem.\ Actually, the Banach--Steinhaus Theorem implies that there is a
dense $G_\delta$-set of such points.\ Being power unbounded, however, is not
enough to ensure that there exists a vector $x_0$ for which
${\lim_n\|T^nx_0\|\kern-1pt=\kern-1pt\infty}.$ (A detailed discussion
along this line can be found in, e.g.,
\cite[Sections III.2 and III.4, pp.48,66]{Bea}$\kern.5pt.)$

\vskip6pt
$\!$For lack of a better name, we call a power unbounded operator $T\kern-1pt$
{\it iterated coercive}\/, or simply {\it coercive}\/ for short, if there is a
vector $x_0$ in the normed space $\X$ such that
${\lim_n\|T^n x_0\|\kern-1pt=\kern-1pt\infty}$; otherwise it will be called
{\it noncoercive}\/.\ In any case, if $T$ is power unbounded, then for every
vector $x_0$ such that ${\sup_n\|T^nx_0\|\kern-1pt=\kern-1pt\infty}$, there
exists a subsequence $\{m_j\}$ of the positive integers such that
${\lim_j\|T^{m_j}x_0\|\kern-1pt=\kern-1pt\infty}$, which will be called a
{\it subsequence of coercivity}\/ of $T$ for the vector $x_0.$

\vskip6pt\noi
{\bf Theorem 4.3$.$}
{\it A weakly quasistable operator\/ $T$ on a Banach space\/ $\X$ is either
\begin{description}
\item{$\kern-18pt$\rm(a)$\kern9pt$}
power bounded 
\vskip1pt\noi
\quad or
\vskip1pt\noi
\item{$\kern-18pt$\rm(b$_1$)$\kern4pt$}
noncoercive power unbounded.\
\end{description}
Moreover, if\/ it is a noncoercive power unbounded operator, then
\begin{description}
\item{$\kern-18pt$\rm(b$_2$)$\kern6pt$}
\hbox{$\kern-2.5pt$for} any vector\/ $\,{x_0\in\X}\,$ such that\/
$\,{\sup_n\|T^nx_0\|=\infty},\,$ every pair of subsequences, one of
coercivity for\/ $x_0$ and the other of weak stability for\/ $x_0$, has a
finite set of common entries.\
\end{description}
}

\proof
First recall that if $\{x_n\}$ is a weakly convergent sequence with entries in
a normed space $\X$, then it is bounded; that is,
$$
x_n\wconv x
\;\;\hbox{for some}\;\,
x\in\X 
\;\;\limply\;\;
{\sup}_n\|x_n\|<\infty.                                           \eqno{(\S)}
$$
Now let $\X$ be a Banach space and let $T$ be an operator on $\X.$
\vskip6pt\noi
(a)
The above implication and an application of the Banach--Steinhaus Theorem
ensure the well-known relation:\
$$
\hbox{weak stability implies power boundedness.}
$$
Furthermore, there also exist power-bounded weakly quasistable operators that
are not weakly stable (see, e.g., Example 4.1).\ This settles the
power-bounded case.\

\vskip6pt\noi
(b)
Suppose $T$ is not power bounded.\ Then there exists a vector $x_0$ in $\X$
for~which $\sup_n\!\|T^nx_0\|=\infty.$ So there is a subsequence $\{m_j\}$ of
the positive integers such~that
$$
{\lim}_j\|T^{m_j}x_0\|=\infty.                                 \eqno{(\dag)}
$$
If $T$ is weakly quasistable, then there also exists a subsequence $\{n_j\}$
such that
$$
T^{n_j}x_0\wconv 0
$$
and hence, by $(\S)$,
\vskip-3pt\noi
$$
{\sup}_j\|T^{n_j}x_0\|<\infty.                                 \eqno{(\ddag)}
$$

\vskip6pt\noi
(b$_1$)
If $T$ is a coercive power unbounded operator, then we can take a vector $x_0$
for which the subsequence $\{m_j\}$ in $(\dag)$ may be replaced by the whole
sequence $\{n\}$ of all positive integers to get
$$
{\lim}_n\|T^nx_0\|=\infty,
$$
and so there is no subse\-quence of the positive integers that does not
satisfy the above limit.\ Hence the~inequality in $(\ddag)$ leads to a
contradiction.\ Consequently, if~$T$ is weakly quasistable and power
unbounded, then it is noncoercive.\

\vskip6pt\noi
(b$_2$)
Next suppose $T$ is a noncoercive power unbounded.\ If the subsequences
$\{m_j\}$ and $\{n_j\}$ in $(\dag)$ and $(\ddag)$ have an infinite set of
common pairwise distinct entries, then there exists a subsequence of the
positive integers made of elements from ${\{m_j\}\cap\{n_j\}}$, which is a
subsequence of both $\{m_j\}$ and $\{n_j\}.$ This subsequence of both
$\{m_j\}$ and $\{n_j\}$ leads to another contradiction according to $(\dag)$
and $(\ddag).$                                                           \qed

\vskip6pt\noi
{\sc Remark 4.4$.$}
{\sf Power unbounded weakly quasistable operators.}
\vskip6pt\noi
As is also widely known, if\/ $T$ is a Banach-space operator, then
\begin{eqnarray*}
r(T)\kern-1pt<1
&\kern-6pt\iff\kern-6pt&
\,T^n\kern-3pt\uconv O
\;\;\limply\;\;
T^n\kern-3pt\sconv O
\;\;\limply\;\;
T^n\kern-3pt\wconv O                                                       \\
&\kern-6pt{\;\,\limply}\kern-6pt&\,\;
{\sup}_n\|T^n\|\kern-1pt<\infty\kern-1pt
\;\;\limply\;\;
r(T)\kern-1pt\le\kern-1pt1,
\end{eqnarray*}
where $r(T)$ stands for the spectral radius of $T.$ In particular, strong
stability implies weak stability, which in turn implies power boundedness:\
\begin{eqnarray*}
{\lim}_n\|T^nx\|=0
\;\;\forall\,x\in\X
&\kern-3pt\limply\kern-3pt&
\forall\,x\in\X,
\,\;{\lim}_n|f(T^nx)|=0
\;\;\forall\,f\in\X^*\!                                                    \\
&\kern-3pt\limply\kern-3pt&
{\sup}_n\|T^n x\|<\infty
\;\;\forall\,x\in\X.{\phantom{\Big|}}
\end{eqnarray*}
Also, if $T$ is a power-bounded operator on a normed space $\X$, then strong
quasistability coincides with strong stability \cite[Proposition 4.2]{KV},
that is,
$$
{\sup}_n\|T^n\|\kern-1pt<\kern-1pt\infty
\,\limply\,
\Big\{{\lim}_n\|T^nx\|\kern-1pt=\kern-.5pt0\;\;\forall\,x\in\kern-1pt\X
\iff
{\liminf}_n\|T^nx\|\kern-1pt=\kern-.5pt0\;\;\forall\,x\in\kern-1pt\X\Big\},
$$
differently from the weak case, where there are power-bounded weakly
quasistable operators that are not weakly stable, as we saw in Example 4.1.\
However, even for the strong case, the equivalence collapses for a power
unbounded operator:\
$$
{\liminf}_n\|T^nx\|\kern-1pt=\kern-1pt0\,\;\forall\,x\kern-1pt\in\kern-1.5pt\X
\;{\kern3pt\not\kern-3pt\limply}\;\,
{\sup}_n\|T^n\|\kern-1pt<\kern-1pt\infty
\quad\&\quad
{\sup}_n\|T^n\|\kern-1pt=\kern-1pt\infty
\,\;{\limply}\;\kern.5pt
T^n\kern-1pt\!\notsconv O.
$$
Indeed, there exists a diagonal operator $T$ on the Hilbert space $\ell_+^2\!$
with
$$
\|T^n\|=(\log\,n)^{_{\scriptstyle\frac{1}{2}}}
\,\;\forall\,n\in\NN 
\quad\;\hbox{and}\;\quad
{\inf}_n\|T^nx\|=0\;\;\forall\,x\in\kern-.5pt\ell_+^2
$$
\cite[Section III.4.A, p.66]{Bea}.\ Since $T$ is power unbounded, it is weakly
unstable and so strongly unstable$.$~Since, for every
${x\kern-.5pt\in\kern-.5pt\ell_+^2}\kern-.5pt$,
${\inf_n\|T^nx\|\kern-1pt=\kern-1pt0}$ implies
${\liminf_n\|T^nx\|\kern-1pt=\kern-1pt0}$ (because
${\|T^{n+1}x\|\kern-1pt\le\|T\|\,\|T^nx\|}$), we get
${\liminf_n|\<T^nx\,;y\>|=0}$ for \hbox{every}
${x,y\kern-1pt\in\kern-.5pt\ell_+^2}$ (by the Schwartz inequality).\ So the
operator $T$ is power unbounded and strongly quasistable and so
weakly quasi\-stable.\ Summing up$:$ (i) for power unbounded operators, strong
quasistability does not imply strong stability, and, moreover,~(ii)
$$
\hbox{weak quasistability does not imply power boundedness.}
$$
\vskip-7pt\noi\vskip5pt

\section{Rajchman Measures and Weak Stability}

The $\sigma$-algebra of all Borel subsets of the real line $\RR$ will be
denoted~by~$\Re.\!$ With~$\Ps$ standing for power set,
${\Re_{[\kern1pt0,1)}\!=\Re\cap\Ps([\kern1pt0,1)\kern.5pt)}$ denotes the
$\sigma$-algebra of all Borel sub\-sets of the interval ${[\kern1pt0,1)}$, and
$\Re_\TT$ stands for the $\sigma$-algebra of all Borel~subsets of the unit
circle $\TT.$ All measures in this paper are positive.\ The terms
``absolutely continuous'', ``continuous'', ``singular'', and ``discrete'' for
a given measure are with respect to normalised Lebesgue measure on
$\kern-.5pt\Re_\TT$ (or on $\Re_{[\kern1pt0,1)})$, unless otherwise~stated.\

\vskip5pt
For each integer ${k\in\ZZ}$, consider the function
${{\bf e}_k\!:\TT\kern-1pt\to\TT}$ given by ${{\bf e}_k(z)=z^k}$ for every
${z\in\TT}.$ For simplicity, write $z^k$ instead of ${\bf e}_k$, as usual, and
let the countable collection of (trigonometric) functions
$\{{\bf e}_k\}_{k\in\ZZ}$ be written simply as $\{z^k\}.$

\vskip6pt
A finite measure ${\mu\!:\Re_\TT\!\to\RR}$ is a {\it Rajchman measure}\/ if
$$
\int_\TT z^k d\mu\to0
\;\;\hbox{as}\;\;
|k|\to\infty.
$$
The term comes from the pioneering works of Rajchman on this class of
\hbox{measures} (e.g., \cite[$\S5$]{Raj}$\kern.5pt$).\ Since Borel measures
over a compact metric space are regular~(see, e.g.,
\cite[Corollary 10.6]{BP}$\kern.5pt$), Rajchman measures are regular because
they are \hbox{finite}~(and so Borel) acting over the compact set $\TT$.\
Take the measure ${\eta\!:\Re_{[\kern1pt0,1)}\!\to\RR}$ induced by the
measure ${\mu\!:\Re_\TT\!\to\RR}$ via a function
${\gamma\!:{[\kern1pt0,1)}\kern-1pt\to\!\TT}$, defined by
${\eta(A)=\mu(\gamma(A)\kern.5pt)}$ for every ${A\in\Re_{[\kern1pt0,1)}}.$
Here the function $\gamma$ is given by $\gamma(\alpha)={e^{2\pi i\alpha}}$ for
every ${\alpha\in[\kern1pt0,1)}$, which is measurable, invertible, with a
measurable inverse.\ The above integral sets the Fourier transform
$\widetilde\eta$ of $\eta$, namely,
$$
\widetilde\eta(k)\,=\int_\TT\kern-1pt z^k d\mu(z)
\,=\int_{[\kern1pt0,1)}\!e^{2\pi ik\alpha}\,d\kern1pt\eta(\alpha)
\;\;\;\hbox{for every}\;\;
k\in\ZZ.
$$
\vskip-2pt

\vskip6pt
For a survey on Rajchman measures, see \cite{Lyo} (also
\cite[Chapter IX]{KL}$\kern.5pt$).\ The~\hbox{basic} properties required here
are listed below (see, e.g., \cite[p.364 and Theorem 3.4]{Lyo}$\kern.5pt$).\
\vskip5pt\noi
$$
\vbox
{\noi
{\narrower\narrower
{\it Every absolutely continuous finite measure is Rajchman}
\vskip0pt\noi\hskip14pt
(in particular, the normalised Lebesgue measure on $\Re_\TT$ is Rajchman).
\vskip0pt\noi
{\it Every Rajchman measure is continuous}\/.
\vskip0pt\noi
{\it There exist singular Rajchman measures}
\vskip0pt\noi\hskip14pt
(and so every singular Rajchman measure is singular-continuous).
\vskip0pt\noi
\vskip0pt}
}
$$
\vskip-1pt

That every absolutely continuous finite measure is Rajchman is readily
\hbox{verified} by the Riemann--Lebesgue Lemma and the Radon--Nikodym
Theorem.\ That \hbox{every} \hbox{Rajchman} measure is continuous is a direct
consequence of Neder's answer \cite[Section 1]{Ned} to a question posed by
Riesz \cite[p.315]{Rie}.\

\vskip6pt
The Cantor--Lebesgue measure $\eta_{\kern1pt[\kern1pt0,1)}$ on
$\Re_{[\kern1pt0,1)}$ is a classical example of a singular-continuous measure
on $\Re_{[\kern1pt0,1)}.$ This is the Borel-Stieltjes measure on
$\Re_{[\kern1pt0,1)}$ generated by the Cantor function associated with the
Cantor set in ${[0,1\kern-1pt]}.$ The \hbox{Cantor--Lebesgue} measure
$\mu_{_{{\scriptstyle\TT}}\!}\kern-1pt$ on $\Re_\TT\kern-1pt$ is obtained from
$\eta_{\kern1pt[\kern1pt0,1)}$ as
${\mu_{_{{\scriptstyle\TT}}\!}(E)
\kern-1pt=\kern-1pt\mu_{[\kern1pt0,1)}\kern-.5pt\big(\gamma^{-1}(E)\big)}$
for ${E\in\Re_\TT}$, where the transformation
${\gamma\!:[\kern1pt0,1)\to\kern-1pt\TT}$ was defined above.\ Thus
$\mu_{_{{\scriptstyle\TT}}\!}$ is a singular-continuous measure on
$\Re_\TT.\kern-1pt$ This, however, is not a Rajchman measure.\ Indeed,~it~can
be verified that
${\int_\TT\kern-1pt z^k d\mu_{_{{\scriptstyle\TT}}\!}
=\!\int_\TT\kern-1pt z^{3k}d\mu_{_{{\scriptstyle\TT}}}\!}$
(i.e.,
${\tilde\eta_{_{\,{\scriptstyle\TT}}\!}(k)
=\tilde\eta_{_{\,{\scriptstyle\TT}}\!}(3k)\kern.5pt}$) for all ${k\in\ZZ}.$
The first example of a singular Rajchman measure is due to Menchoff
\cite[Lemma p.433 and Theorem p.435]{Men}, which is a modification of the
Cantor--Lebesgue \hbox{measure}~on~$\Re_\TT.$

\vskip6pt
Rajchman measures can be equivalently defined as
$$
\int_\TT z^nd\mu\to0
\;\;\hbox{as}\;\;
n\to\infty
$$
with $n$ running over the positive integers.\ Indeed,
${\overline{\int_\TT z^nd\mu}
=\!\int_\TT\overline{z^{\,n}}d\mu
=\!\int_\TT z^{-n}d\mu}$
for every ${n\in\NN}$ whenever ${z\in\TT}$, so that
${\int_\TT\kern-1pt z^nd\mu\kern-1pt\to\kern-1pt0}$ as
${n\kern-1pt\to\kern-1pt\infty}$ implies
${\int_\TT\kern-1pt z^kd\mu\kern-1pt\to\kern-1pt0}$ as
${|k|\kern-1pt\to\kern-1pt\infty}$, and the reverse implication is trivial.

\vskip6pt
Let $\mu$ be a finite measure on $\Re_\TT.$ Consider the Hilbert space
$L^2(\TT,\mu)$ of all~complex functions on $\TT$ that are square-integrable
with respect to $\mu.$ Let ${\varphi\!:\TT\to\TT}$ denote the identity
function (i.e., ${\vphi(z)=z}$ $\mu$\hbox{-}a.e.\ for ${z\in\TT}).$ Take the
multiplication operator $U_{\vphi,\mu}$ on $L^2({\TT,\mu})$ induced by the
identity function $\vphi$, also referred to as the {\it position operator}\/
\cite[p.89]{Hal2}, which is defined by
\vskip6pt\noi
$$
(U_{\vphi,\mu}\psi)(z)=\vphi(z)\psi(z)=z\kern1pt\psi(z)
\quad\;\mu\hbox{-a.e.\ for $z\in\TT$}
\quad\;\hbox{for every $\psi\in L^2(\TT,\mu)$}.
$$
\vskip4pt\noi
So ${(U_{\vphi,\mu}\psi)(z)}={z\psi(z)}$ and
${({U^*_{\vphi,\mu}\psi})(z)=\overline z\psi(z)=z^{-1}\psi(z)}$ for every
${\psi\in L^2(\TT,\mu)}$ and every ${z\in\TT}.\kern-1pt$ Thus
${U^*_{\vphi,\mu}U_{\vphi,\mu}\!=U_{\vphi,\mu}U^*_{\vphi,\mu}\!=I}$, and
${U_{\vphi,\mu}}$ is unitary on $L^2({\TT,\mu}).$ From now on, $U_{\vphi,\mu}$
will always denote such a multi\-plication operator (i.e., the position
operator) on ${L^2(\TT,\mu)}$ for some measure $\mu$ on $\Re_\TT.$ With
${\<\;\,\,;\;\,\>}$~standing for the inner product in $L^2({\TT,\mu})$, we get
for every ${\psi,\phi\in L^2({\TT,\mu})}$ and every~${k\in\ZZ}$,
$$
\<U_{\vphi,\mu}^k\psi\,;\phi\>=\!\int_\TT z^k\psi\kern1pt\overline\phi\,d\mu
\qquad\hbox{so that}\qquad
\<U_{\vphi,\mu}^k 1\,;1\>=\!\int_\TT z^kd\mu,
$$
as $U_{\vphi,\mu}$ is invertible.\ For the second identity, also because the
unit function ${1\!:\!\TT\!\to\kern-1pt\TT}$ (i.e., ${1(z)=1}$
$\mu$\hbox{-}a.e.\ for ${z\in\TT}$) lies in ${L^2(\TT,\mu)}$ since the
measure $\mu$ is finite.

\vskip6pt
Before proving the next proposition, we need the following two well-known
auxiliary results, which will be often required in the sequel.\

\vskip6pt\noi
{\sc Remark 5.1$.$}
{\sf A standard application of the Stone--Weierstrass Theorem.}
\vskip6pt\noi
Let $P(\TT)$ denote the set of all polynomials
${q(\;\;,\;\;)\!:\!\TT\!\times\!\TT\kern-2pt\to\CC}$ in $z$ and $\overline z$
(of the form ${\sum_{k,\ell=0}^N\alpha_{k,\ell}\,z^kz^{-\ell}}$ with
${z^k\kern-1pt,{\overline z}^{\,\ell}\kern-1pt\in\kern-1pt\TT}$
and ${\alpha_{k,\ell}\kern-1pt\in\kern-1pt\CC}.)$ This is a classical density
result:\
\vskip6pt\noi
{\narrower\narrower
$P(\TT)$ is dense in ${(L^p(\TT,\mu),\|\;\;\|_p)}$ for every finite measure
$\mu$ in $\Re_\TT$, \\ and therefore ${(L^p(\TT,\mu),\|\;\;\|_p)}$ is
separable, for every ${p\ge1}.$
\vskip0pt}
\vskip6pt\noi
Indeed, since $\TT$ is {\it compact}\/ in $\CC$, the Stone--Weierstrass
Theorem ensures that $P(\TT)$ is dense in the linear space $C(\TT)$ of all
complex-valued continuous functions on $\TT$ equipped with the
sup-norm$.\kern-1pt$ Thus, since the measure $\mu$ is {\it finite}\/,
$P(\TT)$~is~dense in $C(\TT)$ equipped with the norm-$p$ for any ${p\ge1}.$
Also, since on compact metric spaces Borel measures (in particular, finite
measures) are regular, the set $C(\TT)$ is dense in
${(L^p(\TT,\mu),\|\;\;\|_p)}$ for every ${p\ge1}$ (see, e.g.,
\cite[Theorem~29.14]{Bau}$\kern.5pt).$ By transitivity, $P(\TT)$ is dense in
${(L^p(\TT,\mu),\|\;\;\|_p)}.$ Therefore, as $P(\TT)$ is a linear span of a
countable set, ${(L^p(\TT,\mu),\|\;\;\|_p)}$ is separable.\

\vskip6pt\noi
{\sc Remark 5.2$.$}
{\sf A basic change of variables for complex measurable functions.}
\vskip6pt\noi
Let $\mu$ be a measure on $\Re_\TT$, and let $g$ be a nonnegative
$\Re_\TT$-measurable function on $\TT.$ Take the measure $\eta_g$ on
$\Re_\TT\kern-1pt$ generated by $\mu$ via $g$ defined by
${\eta_g(E)\kern-1pt=\!\int_E g\,d\mu}$ for \hbox{every} $E$ in $\Re_\TT.$
Recall the elementary identity
${\int_E hg\,d\mu\kern-1pt=\!\int_E h\,d\kern1pt\eta_g}$ for every $E$ in
$\Re_\TT$ and every nonnegative $\Re_\TT$-measurable function $h$ on
$\TT.\kern-1pt$ Also, for an arbitrary $\CC$-valued $\Re_\TT$-measurable
function $\zeta$ on $\TT$, consider its natural decomposition into nonnegative
parts given by ${\zeta=h_1-h_2+i(h_3-h_4)}$, where $h_1$, $h_2$, $h_3$, and
$h_4$ are nonnegative $\Re_\TT$-integrable functions on $\TT.\kern-1pt$ So we
get the identity for complex~functions:\
$$
\int_E \zeta\kern1pt g\,d\mu
=\!\int_E \zeta\kern1pt d\kern1pt\eta_g
\quad\hbox{for every}\quad
E\in\Re_\TT.
$$
\vskip-2pt

\vskip6pt
The equivalence between (a) and (e) in Proposition 5.3 below was mentioned in
\cite[p$.$1383]{BM} and used in \cite[Theorem 5]{BM},
\cite[Proposition 3.3]{Kub1}, and \cite[Eq.\ (5.3)]{JJKS}.\

\vskip6pt\noi
{\bf Proposition 5.3$.$}
{\it If\/ $\mu$ is a finite measure on\/ $\Re_\TT\!$ and\/
$U_{\vphi,\mu}\kern-1pt$ is the position operator on\/ $L^2(\TT,\mu)$, then
the following assertions are pairwise equivalent.\
\begin{description}
\item{$\kern-6pt$\rm(a)$\kern2pt$}
${\mu}$ is a Rajchman measure.\
\vskip5pt
\item{$\kern-6pt$\rm(b)$\kern1pt$}
$\int_\TT\kern-1pt z^k g\,d\mu\kern-1pt\to\kern-1pt0$ as\/
${|k|\kern-1pt\to\kern-1pt\infty}$ for every\/
${g\kern-1pt\in\kern-1pt L^p(\TT,\mu)}$, for every\/
${p\kern-1pt\ge\kern-1pt1}.$
\vskip5pt
\item{$\kern-6pt$\rm(c)$\kern2pt$}
$\int_\TT\kern-1pt z^k g\,d\mu\kern-1pt\kern-1pt\to0$ as\/
${|k|\kern-1pt\to\kern-1pt\infty}$ for every\/
${g\kern-1pt\in\kern-1pt L^p(\TT,\mu)}$, for some\/
${p\kern-1pt\ge\kern-1pt1}.$
\vskip5pt
\item{$\kern-6pt$\rm(d)$\kern1pt$}
$\int_\TT\kern-1pt z^k g\,d\mu\kern-1pt\to\kern-1pt0$ as\/
${|k|\kern-1pt\to\kern-1pt\infty}$ for an arbitrary positive\/
${g\in L^\infty(\TT,\mu)}.$
\vskip5pt
\item{$\kern-6pt$\rm(e)$\kern2pt$}
$U_{\vphi,\mu}$ is weakly stable.\
\vskip5pt
\item{$\kern-6pt$\rm(f)$\kern2pt$}
${\<U^n_{\vphi,\mu}\psi\,;\psi\>\kern-1pt\to\kern-1pt0}$ as\/
${n\kern-1pt\to\kern-1pt\infty}$ for an arbitrary\/
${\psi\kern-1pt\in\kern-1pt L^\infty(\TT,\mu)}$ with\/ ${0<|\psi|}.$
\vskip5pt
\item{$\kern-6pt$\rm(g)$\kern1pt$}
${\<U^n_{\vphi,\mu}1\,;1\>\kern-1pt\to\kern-1pt0}$ as\/
${n\kern-1pt\to\kern-1pt\infty}.$
\end{description}

\proof
${\rm(a)\limply(b)}.$
Let ${q=\kern-1pt\sum_{k,\ell=0}^N\alpha_{k,\ell}z^kz^{-\ell}\!}$ be an
arbitrary polynomial in $P(\TT)$ so that
${\int_\TT\!z^n q\,d\mu}$ $\kern-1pt=\kern-1pt{\sum_{k,\ell=0}^N
\alpha_{k,\ell}\int_\TT z^n\,z^kz^{-\ell}d\mu}.$
If $\mu$ is a Rajchman measure, then
$$
\int_\TT z^{n+k-\ell}\,d\mu\to 0
\;\;\;\hbox{as}\;\;
n\to\infty
\;\;\hbox{for each}\;\,
k,\ell\in\ZZ.
$$
Thus ${\int_\TT z^n q\,d\mu\to0}$ for every ${q\in P(\TT)}.$ As we saw in
Remark 5.1, $P(\TT)$ is dense in ${(L^p(\TT,\mu),\|\;\;\|_p)}$ for every
${p\kern-1pt\ge\kern-1pt1}.$ Hence ${\lim_n\int_\TT z^n g\,d\mu=0}$, and since
${\overline{\int_\TT z^k g\,d\mu}}=\!{\int_\TT z^{-k}\,{\overline g}\,d\mu}$,
we get ${\int_\TT z^k g\,d\mu\to0}$ as ${|k|\to\infty}$, for every
${g\in L^p(\TT,\mu)}$, for every~${p\kern-1pt\ge\kern-1pt1}.$

\vskip6pt\noi
${\rm(b)\kern-.5pt\limply\kern-.5pt(c)}$
trivially and
${\rm(c)\kern-.5pt\limply\kern-.5pt(d)}$
since ${L^\infty(\TT,\mu)\kern-1pt\subset\kern-1pt L^p(\TT,\mu)}$ for
${p\kern-.5pt\ge\kern-1pt1}$ as $\mu$
\hbox{is~finite}.\

\vskip6pt\noi
${\rm(d)\limply(a)}.$
Suppose (d) holds.\ Take an arbitrary positive ${g\in L^\infty(\TT,\mu)}$, and
let ${g^{-1}\kern-1pt=\frac{1}{g}}$ denote its reciprocal, which is again a
positive function in ${L^\infty(\TT,\mu)}$, so that both $g$ and $g^{-1}\!$
lie in $L^p(\TT,\mu)$ for every ${p\ge1}$, once $\mu$ is finite.\ Since $g$ is
positive and measurable, take the measure $\eta_g$ on $\Re_\TT$ defined by
${\eta_g(E)\kern-1pt=\!\int_E g\,d\mu}$~for~every~$E$ in $\Re_\TT$ so
that $\eta_g$ is absolutely continuous with respect to $\mu$ (i.e.,
${\eta_g\kern-1pt\ll\kern-1pt\mu}$) and finite as $g$ lies in $L^1(\TT,\mu).$
Also, by Remark 5.2 with ${\zeta\!=z^k}\!$, we get
${\int_E z^kd\kern1pt\eta_g\kern-1pt=\!\int_E z^kg\,d\mu}$~for every
${E\in\Re_\TT}$ and ${k\in\ZZ}.$ Hence, with ${E=\kern-1pt\TT}$, (d) implies
${\int_\TT z^kd\kern1pt\eta_g\!\to0}$ as ${|k|\to\infty}$, and so $n_g$ is
Rajchman.\ Applying Remark 5.2 again, now with
${\zeta\!=\kern-1pt z^kg^{-1}}\!$, it follows~that
$$
{\int_E z^k g^{-1}d\kern1pt\eta_g=\!\int_E z^k g^{-1}g\,d\mu
=\!\int_E z^k d\mu}                                            \eqno{(\star)}
$$
for every $E$ in $\Re_T$ and every $k$ in $\ZZ.$ There are two ways to
conclude the proof.\
\vskip4pt\noi
Conclusion 1.\ As $\eta_g\kern-.5pt$ is Rajchman and
${g^{-1}\kern-1pt\in L^1(\TT,\mu)}$,
${\int_\TT z^k g^{-1}d\kern1pt\eta_g\kern-1pt\to0}$ as ${|k|\to\infty}$
because ${\rm(a)\limply(c)}.$ So $\mu$ is Rajchman by $(\star)$ with
${E=\kern-1pt\TT}.$ Thus (a) holds true.\
\vskip4pt\noi
Conclusion 2.\ For ${k=0}$, $(\star)$ says that
${\mu(E)=\int_E d\mu=\!\int_E g^{-1}d\kern1pt\eta_g}$ for every
${E\in\Re_\TT}$, and so $\mu$ is absolutely continuous with respect to
$\eta_g$ (i.e., ${\mu\kern-1pt\ll\kern-1pt\eta_g}).$ Since $\eta_g\kern-1pt$
is Rajchman, it follows that $\mu$ must be Rajchman.\ (This is a result due
to Milicer-Gru\.zewska (cf.\ \cite[p.365 (2.1)$\kern.5pt$]{Lyo}$\kern.5pt$).\
Thus, again, (a) holds true.\

\vskip6pt\noi
${\rm(a)\limply(e)}.$
Take an arbitrary ${\psi\in L^2(\TT,\mu)}.$ Thus
${|\psi|^2\!\in L^1(\TT,\mu)}.$ If $\mu$ is \hbox{Rajchman}, then
${\lim_n\int_\TT z^n|\psi|^2\,d\mu=0}$ since ${\rm(a)\limply(b)}.$
Equivalently, ${\lim_n\int_\TT U_{\vphi,\mu}^n|\psi|^2\,d\mu=0}$, which means
${\lim_n\<U_{\vphi,\mu}^n\psi\,;\psi\>=0}.$ Since this holds for every
${\psi\in L^2(\TT,\mu)}$, and since ${L^2(\TT,\mu)}$ is a complex
Hilbert~space, it follows that ${\lim_n\<U_{\vphi,\mu}^n\psi\,;\phi\>=0}$ for
every ${\psi,\phi\in L^2(\TT,\mu)}$
(cf.\ Remark 3.6(a,b)$\kern.5pt$); that is, $U_{\vphi,\mu}$ is weakly stable.\

\vskip6pt\noi
${\rm(e)\limply(f)}$
since ${L^\infty(\TT,\mu)\subset L^2(\TT,\mu)}$ because $\mu$ is finite.\

\vskip6pt\noi
${\rm(f)\limply(a)}$
because ${\<U^n_{\vphi,\mu}\psi\,;\psi\>}=\!\int_\TT z^n|\psi|^2 d\mu$ and
${\rm(d)\limply(a)}.$

\vskip6pt\noi
${\rm(g)\!\iff\!(a)}$
as ${\!\int_\TT z^n d\mu=\<U^n_{\vphi,\mu}1\kern1pt;\kern-1pt1\>}$
and ${1\in L^2(\TT,\mu)}$ because $\mu$ is finite.\                      \qed

\section{Rajchman Measures and Trigonometric Basis}

The celebrated Riemann--Lebesgue Lemma says that
$$
\centerline{
${\int_\TT z^k\psi(z)\,d\lambda\to0}\;$ as $\;{|k|\to\infty}\;$ for every
$\,{\psi\in L^2(\TT,\lambda)}$}
$$
\vskip-1pt\noi
$$
\centerline{
if $\lambda$ is the normalised Lebesgue measure on $\Re_\TT$.}
$$
\vskip3pt\noi
There are different proofs of the Riemann--Lebesgue Lemma$.$~One~of them that
uses a bit of operator theory, whose argument is required in the proof of
Theorem 6.1 below, is this.\ As observed in Remark 5.1, since $\TT$ is compact
and $\lambda$ is finite, $L^2(\TT,\lambda)$ is separable.\ Take the position
operator $U_{\vphi,\lambda}\!$ on $L^2(\TT,\lambda)$ so that
${\int_\TT z^k\psi(z)\,d\lambda(z)}={\<U_{\vphi,\lambda}^k\psi\,;1\>}$, where
$\psi$ and the unit function $1$ are in $L^2(\TT,\lambda).$ Recall that (i)
$$
\hbox{$\lambda\kern-1pt$ is the normalised $\kern-1pt$Lebesgue measure
$\kern-.5pt\limply\!$ $\{z^k\}\!$ is an orthonormal basis for
$\kern-1ptL^2(\TT\kern-1pt,\kern-1pt\lambda)$}
$$
(this is a standard result; see, e.g., \cite[p.18]{Hal2}$\kern.5pt$), that
(ii) if an operator shifts an ortho\-normal basis that runs over $\ZZ$, then
it is a bilateral shift of multiplicity one, and also that (iii) every
bilateral shift on a Hilbert space is weakly stable.\
Thus, since $U_{\vphi,\mu}$ shifts $\{z^k\}_{k\in\ZZ}$ independently of the
measure $\mu$~on~$\Re_\TT$,
$$
\hbox{$\{z^k\}_{k\in\ZZ}$ is an orthonormal basis for $L^2(\TT\kern-.5pt,\mu)$
$\;\limply\;$ $U_{\vphi,\mu}$ is weakly stable.}
$$
Therefore, ${\<U_{\vphi,\lambda}^n\psi;\,\phi\>\to0}$ as ${n\to\infty}$ for
every ${\psi,\phi\in L^2(\TT,\lambda)}$, and so
${\int_\TT z^n d\lambda}={\<U_{\vphi,\lambda}^k1;\,1\>\to0}$ as
${|k|\to\infty}.$ This paves the way for the next result.
\vskip2pt

\vskip6pt\noi
{\bf Theorem 6.1$.$}
{\it Let\/ $\mu$ be a measure on}\/ $\Re_\TT.$
\begin{description}
\item{$\kern-9pt$\rm(a)}
{\it If\/ $\{z^k\}$ is an orthonormal basis for\/ $L^2(\TT,\mu)$, then\/
$\mu$\/ is a Rajchman probability measure}\/.\
\end{description}
{\it The converse of\/ {\rm(a)} fails}\/.\
\begin{description}
\item{$\kern-9pt$\rm(b)}
{\it If a probability measure\/ $\mu$ on\/ $\Re_\TT\kern-1pt$ is
singular-continuous and concentrated on the first quadrant of\/
$\kern-.5pt\TT\kern-.5pt$, then\/ $\{z^k\}$ is a set of pairwise nonorthogonal
unit vectors in}\/ ${L^2(\TT,\mu)}.$
\end{description}
{\it This implies the next result}\/.\
\begin{description}
\item{$\kern-7.5pt$\rm(c)}
{\it There exist Rajchman\/ probability\/ measures\/ $\mu$\/ on\/ $\Re_\TT$
for which\/ the set $\{z^k\}$ of unit vectors is not orthogonal in}\/
$L^2(\TT,\mu).$
\end{description}
\vskip-5pt\noi

\vskip5pt
\proof
(a$_1$)
With ${{\bf e}_k(z)=z^k\in\TT}$ for every ${z\in\TT}$ and $\{{\bf e}_k\}$
being an orthonormal set, we get
$1=\|{\bf e}_k\|^2\!
=\!\int_\TT|{\bf e}_k(z)|^2 d\mu(z)
=\!\int_\TT|z^k|^2 d\mu(z)
=\!\int_\TT d\mu
=\mu(\TT)$
for every $k$ in $\ZZ$, so that $\mu$ is a probability measure.\
\vskip6pt\noi
(a$_2$)
Using the preceding argument, we can verify that ${\int_\TT z^n d\mu\to0}$ as
${n\to\infty}.$ In fact, let $U_{\vphi,\mu}\kern-1pt$ be the position operator
on $L^2({\TT,\mu})$.\ If $\{z^k\}$ is an orthonormal basis for
$L^2(\TT,\mu)$, then $U_{\vphi,\mu}$ is a bilateral shift (because it shifts
this orthonormal~basis; i.e., ${U_{\vphi,\mu}{\bf e}_k={\bf e}_{k+1}}$).\
Thus, like every bilateral shift, it is weakly stable, so that
${\<U_{\vphi,\mu}^n 1\,;1\>}={\int_\TT z^n d\mu\to0}$ as ${n\to\infty}.$ This
implies ${\int_\TT z^k d\mu\to0}$ as ${|k|\to\infty}.$ Since $\mu$ is
finite by (a), it is a Rajchman measure.\
\vskip6pt\noi
(b)
The set $\{z^k\}$ is made up of unit vectors because $\mu$ is a probability
\hbox{measure}$.$~Let $\mu$ be singular with respect to the normalised
Lebesgue measure $\lambda$ on $\Re_\TT$, with support
${[\mu]=C\sse\{z\in\TT\!:z=e^{2\pi i\alpha}\;\hbox{for every}\;
\alpha\in[\kern1pt0,\frac{1}{4}]\}}.$
So there exists an $\Re_\TT$-measurable partition $\{B,C\}$ of $\TT$ such that
${\mu(B)=\lambda(C)=0}$ and ${\lambda(\TT)=\lambda(B)}={\mu(C)=\mu(\TT)=1}.$
Hence,
$$
\<z\,;1\>
=\!\int_\TT z\,d\mu(z)
=\!\int_C z\,d\mu(z)
=\!\int_{[\kern1pt0,\frac{1}{4}]}e^{2\pi i\alpha}\,d\kern1pt\eta(\alpha),
$$
with ${\eta\!:\Re_{[\kern1pt0,1)}\to[\kern1pt0,1]}$ given by
${\eta(A)=\mu(\gamma(A)\kern.5pt)}$ for every set ${A\in\Re_{[\kern1pt0,1)}}$,
where ${\gamma\!:[\kern1pt0,1)\to\TT}$ is the homeomorphism (thus a measurable
transformation) of ${[\kern1pt0,1)}$ (equipped with the usual metric in $\RR$)
onto $\TT$ (equipped with the arc-length metric in $\TT$) defined by
${\gamma(\alpha)\kern-1pt=\kern-1pte^{2\pi i\alpha}}$ for every
${\alpha\kern-1pt\in\kern-1pt[\kern1pt0,1)}.$ The support of $\eta$ is
${[\eta]=\gamma^{-1}(\kern.5pt[\mu]\kern.5pt)}$
${=\gamma^{-1}(C)\sse[\kern1pt0,\frac{1}{4}]}.$
If $\mu$ is continuous with respect to $\lambda$, then $\eta$ is continuous
with respect to the restriction to ${[\kern1pt0,1)}$ of the Lebesgue measure
on $\Re$, and therefore singletons in ${[\kern1pt0,\frac{1}{4}]}$ have
$\eta$-measure zero$.$ Therefore, since $\eta$ is concentrated
on~${[\kern1pt0,\frac{1}{4}]}$,
$$
\<z\,;1\>
=\!\int_{(0,\frac{1}{4})}\!\cos(2\pi\alpha)\,d\kern1pt\eta(\alpha)
\,+\;i\!\int_{(0,\frac{1}{4})}\!\sin(2\pi\alpha)\,d\kern1pt\eta(\alpha).
$$
As the functions ${\cos(2\pi\;\;)}$ and ${\sin(2\pi\;\;)}$ are strictly
positive everywhere on ${(0,\frac{1}{4})}$~and
${\eta\big(\kern.5pt(0,\frac{1}{4})\big)
=\eta\big(\kern.5pt[\kern1pt0,\frac{1}{4}]\big)=\mu(C)=1>0}$,
the above integrals are strictly positive$.$~So
$$
\<z\,;1\>\ne0.
$$
A similar argument shows that, for every ${0\ne k\in\ZZ}$,
$$
\<z^k\,;1\>
=\!\int_{(0,\frac{1}{k4})}e^{2k\pi i\alpha}\,d\kern1pt\eta(\alpha)\ne0.
$$
Consider again the position operator $U_{\vphi,\mu}$ on $L^2(\TT,\mu).$ Since
$U_{\vphi,\mu}$ is unitary, and~so is $U_{\vphi,\mu}^j$, we get
${\<z^j\psi\,;z^j\phi\>
=\<U_{\vphi,\mu}^j\psi\,;U_{\vphi,\mu}^j\phi\>
=\<\psi\,;\phi\>}$
for every ${\psi,\phi}$ in ${L^2(\TT,\mu)}$, \hbox{every} $j$ in $\ZZ$, and
every measure $\mu$ on $\Re_\TT.$ Thus, by the above inequality,
$$
\<z^{j+k}\,;z^j\>\ne0
\;\;\hbox{for every}\;\,
j,k\in\ZZ
\;\,\hbox{with}\,\;
k\ne0.
$$
\vskip6pt\noi
(c)
As we saw before, there are singular-continuous Rajchman measures on
$\Re_\TT.$ Thus there are singular-continuous
$\kern-.5pt$Rajchman$\kern-.5pt$ probability$\kern-.5pt$ measures concentrated
on the first quadrant (compress and rotate its support if necessary).\ So
(b) implies~(c).\                                                        \qed

\vskip6pt\noi
{\sc Remark 6.2$.$}
{\sf Just orthogonality is lost.}
\vskip6pt\noi
However, the standard application of the Stone--Weierstrass Theorem summarised
in Remark 5.1 says that, for any finite measure $\mu$ on $\Re_\TT$, functions
on $L^p(\TT,\mu)$~are approximated by the trigonometric polynomials in
$P(\TT).$ In other words, using the argument in the proof of Proposition 5.3
(for ${\rm(a)\limply(b)}$, and so for ${\rm(a)\limply(e)}\kern.5pt)$, for every
${p\kern-1pt\ge\kern-1pt1}$, ${\bigvee_{k\in\ZZ}\!\{z^k\}\!=L^p(\TT,\mu)}$,
where $\bigvee_{k\in\ZZ}\!\{z^k\}$ is the closure of the linear span of
$\{z^k\}$ (see, e.g., \cite[Theorem 7.2]{Con} for the case of Lebesgue
measure~on~$\Re_\TT$).\

\section{Rajchman Measures and Weak Quasistability}

We say that a finite measure $\mu$ on $\Re_\TT\kern-1pt$ is
{\it quasi-Rajchman}\/ if
$$
{\liminf}_{|k|\to\infty}\Big|\!\int_\TT\!z^k d\mu\Big|=0.
\quad\hbox{Equivalently,}\quad
{\liminf}_n\Big|\!\int_\TT\!z^n d\mu\Big|=0,
$$
which means that there exists a subsequence $\kern-1pt\{n_j\}\kern-1pt$ of the
positive integers such~that
$$
\int_\TT z^{n_j}d\mu\to0
\;\;\hbox{as}\;\;
j\to\infty.
$$

\vskip6pt
Let $\mu$ be a finite measure on $\Re_\TT.\kern-1pt$ The Lebesgue
Decomposition Theorem says~that
$$
\mu=\mu_a\!+\mu_s,
$$
where $\mu_a$ and $\mu_s$ are the absolutely continuous and the singular parts
of $\mu$, respec\-tively.\ Moreover, $\mu_s$ can be further decomposed as
${\mu_s=\mu_{sc}\!+\mu_{sd}}$ yielding a refinement of the above decomposition
(see, e.g., \cite[Corollary 7.14]{EMT}$\kern.5pt$),
$$
\mu=\mu_a\!+\mu_{sc}\!+\mu_{sd},
$$
where $\mu_{sc}$ and $\mu_{sd}$ are the singular-continuous and
singular-discrete (or simply discrete) parts of $\mu.$ This is the
{\it canonical decomposition}\/ of $\mu$, where any of theses parts may be
absent.\ The continuous part $\mu_c$ of $\mu$ consists of the absolutely
continuous and singular-continuous parts of $\mu$:\
$$
\mu_c=\mu_a\!+\mu_{sc}.
$$
\vskip-2pt

\vskip6pt
Recall that every Rajchman~measure is continuous.\ Thus, if a measure on
$\Re_\TT$~has a singular-discrete part, then it is not Rajchman.\ That is,
$$
\hbox{every singular-discrete measure on $\Re_\TT$ is not Rajchman.}
$$
In other words, for a singular-discrete measure $\mu_{sd}$ on $\Re_\TT$, we
have ${\int_\TT z^n d\mu_{sd}\not\to0}$ as ${n\to\infty}.$ Then there exists
at least one subsequence $\{m_j\}$ of the positive integers for which
${\int_\TT z^{m_j}d\mu_{sd}\not\to0}$ as ${j\to\infty}.$ This is not the case
for every subsequence$.$~In~fact,
$$
\hbox{singular-discrete measures on $\Re_\TT$ may be quasi-Rajchman or not.}
$$
The simplest examples are these.\ Set ${\mu_{sd}=\delta_1}$, the Dirac measure
at ${1\in\TT}$, so that ${\int_\TT z^n d\mu_{sd}=z^n|_1=1}$ for all
${n\in\NN}$, and hence there is no subsequence $\{n_j\}$ for which
${\lim_j\kern-1pt\int_\TT z^{n_j}d\mu_{sd}\kern-1pt=0}.$ On the other end, set
${\mu_{sd}=\delta_1+\delta_{-1}}$ so that
${\int_\TT z^n d\mu_{sd}}={z^n|_1+z^n|_{-1}=1+(-1)^n}$ for every ${n\in\NN}$,
and hence if $\{n_j\}$ is the subsequence of all odd integers, then
${\int_\TT z^{n_j}d\mu_{sd}=0}$ for every ${j\in\NN}.$

\vskip6pt
$\!$Again, every Rajchman measure is continuous, but the converse fails.\
As we saw in Section 5, the Cantor--Lebesgue measure on $\Re_\TT$ is
singular-continuous but not Rajchman.\ A continuous measure, however, is
always quasi-Rajchman.

\vskip6pt\noi
{\bf Lemma 7.1$.$}
{\it If a finite measure on\/ $\Re_\TT$ is continuous, then it is
quasi-Rajchman}\/.\

\proof
Let $\mu$ be a finite measure on $\Re_\TT.$ The Wiener characterisation
\cite[Eq.(31) p.81]{Wie} of continuous measures (see also, e.g.,
\cite[p.365]{Lyo}$\kern.5pt$) says that $\mu$ is continuous if and only if
the following mean of bounded summands goes to zero:
\goodbreak\vskip4pt\noi
$$
\hbox{$\mu$ is continuous}
\;\;\iff\;\;
\smallfrac{1}{2k+1}{\sum}_{|\ell|\le k}\Big|\int_\TT\!z^\ell d\mu\Big|\to0
\;\;\hbox{as}\;\;
|k|\to\infty,
$$
\vskip2pt\noi
that is,
${\smallfrac{1}{2n+1}\Big(\Big|\!\int_\TT d\mu\Big|
\kern-1pt+\kern-1pt\sum_{\ell=1}^n\kern-1pt\Big|\!\int_\TT z^\ell d\mu\Big|
\kern-1pt+\kern-1pt\sum_{\ell=1}^n\kern-1pt\Big|\!\int_\TT{\overline z}^\ell
d\mu\Big|\Big)\kern-1pt\to\kern-1pt0}$
as ${n\kern-1pt\to\kern-1pt\infty}.$~\hbox{Therefore},
$$
\hbox{$\mu$ is continuous}
\;\;\iff\;\;
\smallfrac{1}{2n+1}
{\sum}_{\ell=1}^n\Big|\int_\TT z^\ell\!d\mu\Big|\to0
\;\;\hbox{as}\;\;
n\to\infty,
$$
since ${\big|\!\int_\TT d\mu\big|
\kern-1pt=\kern-1pt\!\int_\TT d\mu
\kern-1pt=\kern-1pt\mu(\TT)}.$
If $\mu$ is not quasi-Rajchman, then
${\liminf_n\big|\!\int_\TT\!z^n d\mu\big|\kern-1pt>\kern-1pt0}$, which means
that there is an ${\veps\kern-1pt>\kern-1pt0}$ and a positive integer
$\ell_\veps$ such that
${\big|\!\int_\TT z^{\ell}d\mu\big|\kern-1pt>\kern-1pt\veps}$
for every ${\ell\kern-1pt>\kern-1pt\ell_\veps}.$
Hence, for ${n>\ell_\veps}$,
\begin{eqnarray*}
\smallfrac{1}{2n+1}{\sum}_{\ell=1}^n\Big|\int_\TT z^\ell\!d\mu\Big|
&\kern-6pt=\kern-6pt&
\smallfrac{1}{2n+1}
\Big({\sum}_{\ell=1}^{\ell_\veps}\Big|\int_\TT z^\ell\!d\mu\Big|
+{\sum}_{\ell=\ell_\veps+1}^n\Big|\int_\TT z^\ell\!d\mu\Big|\Big)          \\
&\kern-6pt>\kern-6pt&
\smallfrac{1}{2n+1}
\Big({\sum}_{\ell=1}^{\ell_\veps}\Big|\int_\TT z^\ell\!d\mu\Big|\Big)
+\smallfrac{n+\ell_\veps}{2n+1}\veps
\ge\smallfrac{n+\ell_\veps}{2n+1}\veps
>\smallfrac{1}{3}\veps,
\end{eqnarray*}
so that $\mu$ is not continuous. Thus if $\mu$ is continuous, then it is
quasi-Rajchman.                                                          \qed

\vskip6pt
The Cantor--Lebesgue measure on $\Re_\TT$, which is singular-continuous and
was shown in Section 5 not to be Rajchman, is quasi-Rajchman by Lemma~7.1.\

\vskip6pt
Since a measure is continuous if and only if it has no singular-discrete part
in its canonical decomposition, an immediate consequence of Lemma 7.1 reads
as follows.

\vskip6pt\noi
{\bf Corollary 7.2$.$}
{\it If a finite measure $\mu$ on\/ $\Re_\TT$ has no singular-discrete part,
then it is quasi-Rajchman}\/.\

\vskip6pt
The converse fails$:$ as we have seen before, there are singular-discrete
finite meas\-ures $\eta$ on $\Re_\TT$ for which
${\lim_j\!\int_\TT z^{n_j}d\,\eta=0}$ for some subsequence $\{n_j\}.$

\vskip6pt
Proposition 5.3 shows that $U_{\vphi,\mu}$ is weakly stable if and only if the
measure $\mu$ on $\Re_\TT$ is a Rajchman measure.\ Thus, if $\mu$ is
absolutely continuous, then $U_{\vphi,\mu}$ is weakly stable.\ The next result
establishes a counterpart for continuous measures.

\vskip6pt\noi
{\bf Theorem 7.3$.$}
{\it If a finite measure\/ $\mu$ on\/ $\Re_\TT\!$ is continuous, then the
position operator\/ $U_{\vphi,\mu}\kern-1pt$\/ on $L^2(\TT,\mu)$ is weakly
quasistable}\/.

\proof
Let $\mu$ be a finite continuous measure on $\Re_\TT\kern-1pt$, take a
measurable nonnegative function $g$ on $\TT$, and consider the measure
$\eta_g\kern-1pt$ on $\Re_\TT\kern-1pt$ given by
${\eta_g(E)\kern-1pt=\!\int_E g\,d\mu}$ for
${E\kern-1pt\in\kern-1pt\Re_\TT}.\kern-1pt$ Thus $\eta_g$ is absolutely
continuous with respect to $\mu.$ So $\eta_g\kern-1pt$ is continuous with
respect to Lebesgue measure on $\Re_\TT$ because $\mu$ is.\ Also, suppose the
nonnegative function $g\kern-1pt$ lies in ${L^1(\TT,\mu)}.$ Then $\eta_g$ is
finite.\ Hence, according to Lemma 7.1,
$$
\hbox{$\eta_g$ is quasi-Rajchman.}                                 \eqno{(*)}
$$
Moreover, for an arbitrary $\CC$-valued measurable function $\zeta$, we get
(cf.\ Remark~5.2)
$$
\int_\TT \zeta\kern1pt g\,d\mu
=\!\int_\TT \zeta\kern1pt d\kern1pt\eta_g.                        \eqno{(**)}
$$
By $(*)$, there is a subsequence $\{n_j\}$ for which
${\lim_j\int_\TT z^{n_j}d\kern1pt\eta_g=0}$, and so by $(**)$
$$
\int_\TT z^{n_j}g\,d\mu
=\!\int_\TT z^{n_j}d\kern1pt\eta_g\to0
\;\;\hbox{as}\;\;
j\to\infty.
$$
Since this holds for an arbitrary nonnegative function $g$ in $L^1(\TT,\mu)$,
take an arbitrary function $\psi$ in $L^2(\TT,\mu)$ so that $|\psi|^2$ is a
nonnegative function in $L^1(\TT,\mu).$~Then
\goodbreak\vskip4pt\noi
$$
\<U_{\vphi,\mu}^{n_j}\psi\,;\psi\>
=\<z^{n_j}\psi\,;\psi\>
=\!\int_\TT z^{n_j}|\psi|^2\,d\mu\to0
\;\;\hbox{as}\;\;
j\to\infty\quad\hbox{for every}\quad
\psi\in L^2(\TT,\mu).
$$
\vskip2pt\noi
Since $L^2(\TT,\mu)$ is a complex Hilbert space, Remark 3.5 ensures that
$$
\<U_{\vphi,\mu}^{n_j}\psi\,;\phi\>\to0
\;\;\hbox{as}\;\;j\to\infty
\quad\hbox{for every}\quad
\psi,\phi\in L^2(\TT,\mu).                                        \eqno{\qed}
$$

\vskip6pt\noi
{\sc Remark $\kern-.5pt7\kern-.5pt.4.$}
{\sf A weakly quasistable operator with a common weak stability
\hbox{subsequence}.}
\vskip6pt\noi
A weakly unstable operator that is weakly quasistable was exhibited in
\hbox{Example}~4.1.\ Theorem 7.3 ensures the existence of a special class of
such operators, namely,
$$
\hbox{weakly unstable unitary operators that are homogeneously weakly
quasistable}.
$$
Indeed, for the position operator $U_{\vphi,\mu}$ considered in the proof of
Theorem 7.3, weak quasi\-stability is achieved with a subsequence $\{n_j\}$
that depends only on the fact that $\mu$ is continuous according to
Lemma 7.1.\ Also, the proof of Lemma 7.1 does~not specify a subsequence
$\{n_j\}$ ensuring that $\mu$ is quasi-Rajchman.\ Thus,~that~subsequence
$\{n_j\}$~of weak stability for $U_{\vphi,\mu}$ in the proof of Theorem 7.3
does not depend on the \hbox{function} $g$, and so it does not depend on the
function $\psi.$ Then the same~$\{n_j\}$ is a subsequence of weak stability of
all $\psi.$ This holds in contrast to the general case of Remark 3.7(b).\ In
other words, and according to Definition 3.3,

\vskip4pt\noi
{\narrower\narrower
the position operator $U_{\vphi,\mu}$ for a finite non-Rajchman continuous
measure $\mu$ (in particular, if $\mu$ is not absolutely continuous) is
unitary (so power bounded) and homogeneously weakly quasistable but not weakly
\hbox{stable}.
\vskip1pt}

\vskip6pt
We close the paper by giving a characterisation for non-Rajchman measures or,
equivalently, for weakly unstable position operators.\

\vskip6pt\noi
{\bf Corollary 7.5$.$}
{\it Let\/ $\mu$ be a finite continuous measure on\/ $\Re_\TT.$ If\/ $\mu$ is
not \hbox{Rajchman}, then there exists a nontrivial subsequence\/ $\{m_j\}$ of
the positive integers and a subsequence\/ $\{n_j\}$ of\/ $\{m_j\}$ such that 
$$
\int_\TT\!z^{m_j}\,d\mu\not\to0
\;\;\hbox{as}\;\;
j\to\infty
\qquad\!\hbox{\it but}\qquad
{\lim}_{\,j}\!\int_\TT\!z^{n_j}\,d\mu=0.
$$
\vskip-2pt\noi

\proof
If a finite measure on $\Re_\TT$ is not Rajchman, then
${\int_\TT z^nd\mu\not\to0}$ as ${n\to\infty}.$ So there is a subsequence
$\{{m'\kern-2pt}_j\}$ of the positive integers for which
$$
\int_\TT\!z^{{m'\kern-2pt}_j}\,d\mu\not\to0
\;\;\;\hbox{as}\;\;\;
j\to\infty.
$$
According to Lemma 7.1, $\mu$ is quasi-Rajchman.\ So there is a subsequence
$\{n_j\}$ of the positive integers for which $\{n_j\}$ and all its
subsequences are such that
$$
\int_\TT\!z^{n_j}\,d\mu\to0
\;\;\;\hbox{as}\;\;\;
j\to\infty.
$$
\vskip6pt\noi
(a)
Now identify the subsequences $\{{m'\kern-2pt}_j\}$ and $\{n_j\}$ of the
positive integers with subsets of $\NN.$ Consider the subset
${\{{m'\kern-2pt}_j\}\!\cup\!\{n_j\}}$ of $\NN$ made up of all elements from
$\{{m'\kern-2pt}_j\}$ and $\{n_j\}.$ As a union of sets, the possibly common
elements from $\{{m'\kern-2pt}_j\}$ and $\{n_j\}$ appear only once in
${\{{m'\kern-2pt}_j\}\!\cup\!\{n_j\}}.$ Let this set be properly ordered so as
to make it strictly increasing, and identify it with a subsequence of the
positive integers, say $\{m_j\}.$ Thus $\{m_j\}$ is a supersequence of both
$\{{m'\kern-2pt}_j\}$ and $\{n_j\}$ and so
$$
\int_\TT\!z^{m_j}\,d\mu\not\to0
\;\;\;\hbox{as}\;\;\;
j\to\infty.
$$
(b)
If ${\{{m'\kern-2pt}_j\}\!\cup\!\{n_j\}}$ is nontrivial, we are done.\ If not,
we can make it nontrivial as fol\-lows.\ First recall that an arbitrary
subsequence of the positive integers is non\-trivial if and only if its
complement in $\NN$ is another (infinite) subsequence.\ Thus the sub\-sequence
$\{n_j\}$ is certainly nontrivial; otherwise, $\mu$ would be Rajchman.\ Thus
there exists~a pair of complementary subsequences of $\{n_j\}$, say,
$\{{n'\kern-2pt}_j\}$ and $\{{n''\kern-3pt}_j\}$, both non\-trivial with
${\{{n'\kern-2pt}_j\}\cup\{{n''\kern-3pt}_j\}=\{n_j\}}$
and ${\{{n'\kern-2pt}_j\}\cap\{{n''\kern-3pt}_j\}=\void}$, such that
${\int_\TT\!z^{{n'\kern-2pt}_j}\,d\mu\to0}$ and
${\int_\TT\!z^{{n''\kern-3pt}_j}\,d\mu\to0}$ as ${j\to\infty}.$ Then
 ${\{m'_j\}\!\cup\!\{n'_j\}}$ is nontrivial, and the
argument in item (a) holds for $\{n_j\}$ replaced with $\{n'_j\}.$       \qed

\vskip6pt
Nontriviality for the sequence $\{m_j\}$ and the fact that $\{n_j\}$ is a
subsequence of it are crucial in the above result, as well as in the next
one.\ Corollary 7.5 can be rephrased in terms of the position operator
${U_{\vphi,\mu}}$ on $L^2(\TT,\mu)$ according to Proposition 5.3(a,e) and
Theorem 7.3, as follows.\

\vskip6pt\noi
{\bf Corollary 7.6$.$}
{\it Let\/ $\mu$ be a finite continuous measure on\/ $\Re_\TT.$ If\/ the
position operator\/ ${U_{\vphi,\mu}}$\/ on $L^2(\TT,\mu)$ is weakly unstable,
then there exists a nontrivial subsequence\/ $\{m_j\}$ of the positive
integers and a subsequence\/ $\{n_j\}$ of\/ $\{m_j\}$ such that
$$
U^{m_j}_{\vphi,\mu}\psi\kern-1pt\notwconv0
\;\;\hbox{for some}\;\;
\psi\in L^2(\TT,\mu)
\quad\;\hbox{\it but}\;\quad
U^{n_j}_{\vphi,\mu}\psi\kern-1pt\wconv0
\;\;\hbox{for every}\;\;
\psi\in L^2(\TT,\mu).
$$
}

\vskip-10pt\noi
\section*{Acknowledgement}

I thank my co-authors in reference \cite{JJKS} for several discussions
we had on Rajchman measures and weak stability.\

\vskip-10pt\noi
\bibliographystyle{amsplain}

\end{document}